\documentclass{amsart}
\usepackage[utf8]{inputenc}
\usepackage[T1]{fontenc}
\usepackage[british]{babel}
\usepackage{lmodern,microtype}
\usepackage[a4paper,margin=2.5cm]{geometry}
\usepackage{amsmath,amssymb,amsthm,mathtools,bm}
\usepackage{graphicx,booktabs,array,tabularx,adjustbox}
\usepackage{caption,subcaption,enumitem}
\usepackage[hidelinks]{hyperref}
\usepackage[capitalise,nameinlink]{cleveref}
\setlist{itemsep=2pt,topsep=4pt}
\allowdisplaybreaks
\numberwithin{equation}{section}
\newtheorem{proposition}{Proposition}[section]

\theoremstyle{remark}
\newtheorem{remark}[proposition]{Remark}
\newcommand{\R}{\mathbb{R}}
\newcommand{\C}{\mathbb{C}}
\newcommand{\Z}{\mathbb{Z}}
\newcommand{\E}{\mathbb{E}}
\newcommand{\norm}[1]{\lVert #1\rVert}
\newcommand{\radv}{\widehat r_{\mathrm{adv}}}
\newcommand{\Nadv}{\widehat N_{\mathrm{adv}}}

\newcommand{\ghlogo}{\raisebox{-0.2em}{\includegraphics[height=1em]{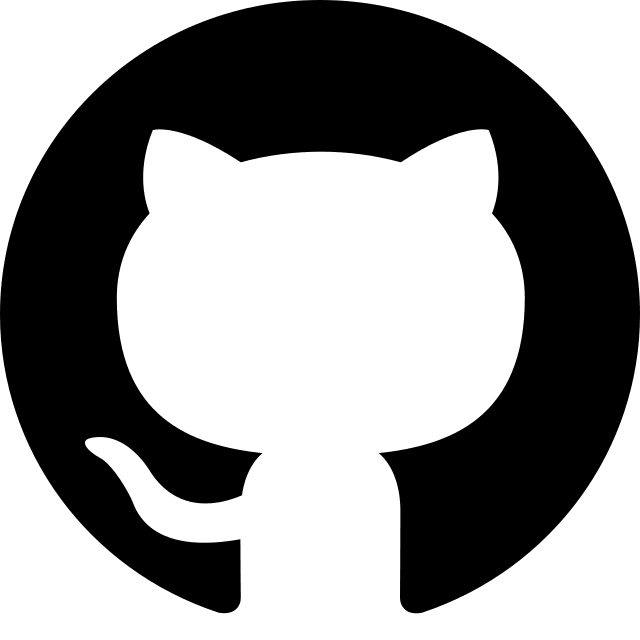}}}
\DeclareMathOperator{\Pf}{Pf}
\DeclareMathOperator{\diag}{diag}
\DeclareMathOperator{\relu}{relu}
\title{\bfseries Learning Almost Complex and Complex Structures}

\author{Edward Hirst}
\author{Tancredi Schettini Gherardini}
\author{Tom\'as S. R. Silva}
\author{Jakob Stein}

\address[E. Hirst, T. S. R. Silva, J. Stein]{%
Instituto de Matem\'atica, Estat\'istica e Computa\c{c}\~ao
Cient\'ifica (IMECC), Universidade Estadual de Campinas
(UNICAMP), Campinas, Brazil}

\email[E. Hirst]{%
\href{mailto:ehirst@unicamp.br}{ehirst@unicamp.br}}
\email[T. S. R. Silva]{%
\href{mailto:tomas@ime.unicamp.br}{tomas@ime.unicamp.br}}
\email[J. Stein]{%
\href{mailto:jstein@unicamp.br}{jstein@unicamp.br}}

\address[T. Schettini Gherardini]{%
Mathematisches Institut, University of Bonn, Bonn, Germany}
\address[T. Schettini Gherardini]{%
Max Planck Institute for Mathematics, Bonn, Germany}

\email[T. Schettini Gherardini]{%
\href{mailto:tsg@math.uni-bonn.de}{tsg@math.uni-bonn.de}}

\date{\today}
\begin{document}
\hypersetup{pdftitle={Learning Almost Complex and Complex Structures},
  pdfauthor={Edward Hirst; Tancredi Schettini Gherardini; Tomas S. R. Silva; Jakob Stein},
  pdfsubject={Numerical searches and geometric reconstructions on spheres and products of spheres}}

\begin{abstract}
We study the numerical constructions of almost complex and complex structures on spheres and products of spheres up to real dimension six. Our general set-up parametrises almost complex structures via a global matrix field in an ambient Euclidean space, compressed by the tangent projector. The projector's algebraic residual measures the almost complex condition, while its Nijenhuis tensor is evaluated in the intrinsic norm of the round product metric. 
Polynomial features and neural networks are trained without prescribed targets to perform general searches for these structures. The almost complex search passes numerical admissibility checks on all nine examples that admit such structures and rejects $S^4$, where a topological obstruction is known. 
A matrix correction then enforces the almost complex identity, subject to a novel uniform admissibility condition. For integrability, a free polynomial representation with Levenberg--Marquardt optimisation learns the complex structures on $S^2\times S^2$, $S^3\times S^1$, $S^3\times S^3$ and $S^5\times S^1$, with the largest Nijenhuis norms found by adversarial search all below $10^{-9}$. These general searches were unsuccessful on $S^6$, $S^4\times S^2$ and $S^4\times T^2$. Therefore, for these three cases direct numerical reconstructions were built in torus-fibred coordinates. They are integrable by construction, satisfy their compatibility conditions to better than $10^{-11}$, and corroborate the recent constructions on $S^6$ and $S^4\times S^2$.\\
\small Code: \ghlogo \, \url{https://github.com/edhirst/learning-acs-cs}
\end{abstract}
\maketitle
\thispagestyle{empty}
\clearpage
\setcounter{page}{1}
\hypersetup{pageanchor=true}
\setcounter{tocdepth}{2}
\tableofcontents
\clearpage

\section{Introduction}\label{sec:intro}
\subsection{Complex structures on spheres and their products}
A natural setting to search for complex manifolds realising given smooth structures is on spheres, and their products. For instance, Hopf manifolds, Calabi--Eckmann manifolds, and the product of Riemann spheres are all classical examples of complex manifolds with a sphere product as the underlying space~\cite{Hopf48,CE53,Samelson53}. Additionally, this area is especially topical due to a recent spate of preprints which claim novel complex structure constructions on $S^6$ ~\cite{AC26,Engel26}, and $S^4 \times S^2$ ~\cite{GFL26,LLZ26}. 

This work initiates a programme for numerical search of almost complex and complex structures on general manifolds, here with focus on spheres and products of spheres.
This numerical programme uses and compares both polynomial and neural network ansatze, testing up to real dimension six (see Table \ref{tab:manifolds}); with the intention that the approach extends naturally to higher dimensions and other manifolds.

The sphere products focused on in this work were specifically those with two factors, however at the time of finishing the numerical work the existence of a complex structure on $S^4 \times T^2$ was unknown, and hence added to the list. This was the final remaining open candidate in dimension up to six for all products of spheres. During our numerical investigation, we also studied an explicit period model for a complex structure on $S^4 \times T^2$, detailed in \cref{sec:s4t2-periods}; a complex structure on this manifold has since been constructed in \cite{Gao:2026tyy}, giving an affirmative answer to this genus-one case of Calabi's problem. Future work will include a detailed derivation and a computer-assisted verification of a complex structure on this manifold.

\subsection{Mathematical formulation} Our primary strategy in our direct search is to view a complex structure on a smooth manifold $M$ as a tensor $J$ that satisfies $J^2 = -I$, i.e. $J$ is an \textit{almost complex structure} (ACS), such that the associated Nijenhuis tensor $N_J$ vanishes, cf. \S \ref{sec:background} equation \eqref{eq:nijenhuis}. By fixing an embedding of $M$ into some suitable $n$-dimensional Euclidean space, we can describe any endomorphism $J$ using a smooth map from $M$ to the space of $n$ by $n$ matrices\footnote{Although this description is not unique, this did  not pose a significant problem for our numerical search.}. The standard Euclidean norm on matrices then descends to a well-defined pointwise tensor norm on $M$ via the projection, and our goal becomes minimising $J^2 + I$ and $N_J$ with respect to these norms.  

\begin{remark} Previous iterations of this project used an atlas of charts to describe the manifolds, as in~\cite{Hirst:2025seh}, and $GL(n,\C)$-structures to describe the almost complex structures, but this had worse numerical convergence, took longer to run, and had less geometric features baked into the architecture (see \cite{gherardini2026minimalsurfacesknotsneural, Hirst:2026cpm, users_guide} for a discussion on these themes). The current approach is both simpler and more easily adaptable to arbitrary manifolds.    
\end{remark}

\begin{table}[h]
\centering\small
\caption{The manifolds studied in this work. Here $d$ is the real dimension and $N$ the ambient dimension of the product embedding. ACS and CS record whether the manifold admits an almost complex and a complex structure respectively; $^*$ marks complex structures claimed in recent preprints.}
\label{tab:manifolds}
\begin{tabularx}{\textwidth}{@{}lccccX@{}}
\toprule
Manifold & $d$ & $N$ & ACS & CS & Notes\\
\midrule
$S^2$ & 2 & 3 & yes & yes & Riemann sphere; surface control\\
$T^2$ & 2 & 4 & yes & yes & elliptic curve; surface control\\
$S^4$ & 4 & 5 & no & no & \cite{BS53}; no ACS; topological obstruction\\
$S^2\times S^2$ & 4 & 6 & yes & yes & $\mathbb{CP}^1\times\mathbb{CP}^1$; integrability control\\
$S^3\times S^1$ & 4 & 6 & yes & yes & \cite{Hopf48}; Hopf surface; integrability control\\
$S^6$ & 6 & 7 & yes & yes$^*$ & \cite{AC26,Engel26}; recent preprints; difficult ambient search\\
$S^3\times S^3$ & 6 & 8 & yes & yes & \cite{CE53,Samelson53}; Calabi--Eckmann; integrability control\\
$S^5\times S^1$ & 6 & 8 & yes & yes & \cite{Hopf48}; Hopf manifold; integrability control\\
$S^4\times S^2$ & 6 & 8 & yes & yes$^*$ & \cite{GFL26,LLZ26}; recent preprints; difficult ambient search\\
$S^4\times T^2$ & 6 & 9 & yes & yes$^*$ & \cite{Gao:2026tyy}; recent preprint; symmetry and reconstruction\\
\bottomrule
\end{tabularx}
\end{table}

\subsection{Numerical approaches and neural networks} The numerical approach is motivated by recent progress in learning geometric structures from their defining equations. Neural representations have been used to approximate Calabi--Yau metrics~\cite{Douglas:2020hpv, Larfors:2022nep}, Einstein metrics~\cite{Hirst:2025seh, Cortes:2026kfx, SchettiniGherardini:2026bdb}, and $\mathrm{G}_2$ metrics~\cite{Douglas:2024pmn, Heyes:2026rch}; including examples where explicit formulae are unavailable. In the physics-informed approach~\cite{Raissi19}, a highly-parametrised function ansatz is iteratively trained to identify and approximate a solution to a PDE by evaluating differential residuals at sampled points. 
The learning process is deliberately completely general, the series of initial ansatze include polynomials of increasing order and neural networks where the latter are universal function approximators and hence with enough computational resources should be able to find solutions where they exist with no prescribed structure or a prior knowledge. This offers a flexible way to search for $J$ across general manifolds, including in cases where existence is unknown. 

The central results set up these ansatze for each considered manifold, such that they take as input the coordinates of a point on the manifold in the coordinate system of the standard embedding ($S^n \subset \mathbb{R}^{n+1}$ for each sphere). 
Training of the parameters in these ansatz involved iterative updating of their values to minimise the almost complex residual in the first stage, then after an explicit correction to the ansatz another stage of training is run to minimise the Nijenhuis residual.
Throughout the unsupervised searches no information was provided about literature-known structures; these analytic structures already known in the literature were only used to independently check the differential operators and the expressive power of the ansatz function classes. 
Where these general learning methods struggled, specifically for these novel delicate six dimensional structures, a more specific ansatz was implemented for the numerical approximation using information about the fibration, monodromy, and functional equations of a geometric construction; where then a spectral solver learnt the compatible period data. 
Both approaches minimise residual equations, treated as losses in the traditional AI methodology, rather than fit explicit labelled solutions, but the second incorporates substantially more geometric information. We therefore distinguish general searches from geometry-guided reconstructions throughout.

All numerical methods implemented represent $J$ globally. If $P(p)$ is the orthogonal projector onto the tangent space of the embedded manifold, a matrix-valued function $M_\theta(p)$ then defines a tangential field $K_\theta=P M_\theta P$, a genuine endomorphism of the tangent bundle and hence a candidate $J$. The projected ambient coordinate vectors form a Parseval frame, allowing the Nijenhuis residual to be evaluated with the intrinsic norm of the round product metric. For the initial almost complex stage, a sequence of increasingly expressive polynomial classes is combined with neural refinement. A uniform algebraic residual below a fixed threshold would guarantee that a matrix-function correction produces an exact almost complex structure like the one used for the start of the integrability learning stage; dense sampling and adversarial maximisation provide numerical tests of this condition.

The resulting searches are recorded in \S \ref{sec:results}. We recover almost complex candidates on all nine admissible manifolds in the ten-case study, and cannot recover one for $S^4$, which is known not to admit almost complex structures, see \Cref{tab:acs-paired}. For integrability, we refer to \Cref{tab:classical-learning}, where direct polynomial optimisation with Levenberg--Marquardt reaches an estimated maximum Nijenhuis norm of order $10^{-10}$ on $S^2\times S^2$, $S^3\times S^1$, $S^3\times S^3$ and $S^5\times S^1$. These results are on-par with the outcomes for the two-dimensional examples where, as explained in \S \ref{sec:background}, their $J$'s are automatically integrable once the almost complex condition is satisfied. Neural conjugation provides another way to preserve that condition, but converges less accurately within the recorded budgets; highlighting the benefit of a complementary choice of representation.

The direct searches on $S^6$, $S^4\times S^2$, and $S^4\times T^2$ do not reach the integrability tolerance: their best mean Nijenhuis norms are $0.09$ ($S^4\times T^2$), $0.37$ ($S^6$) and $0.86$ ($S^4\times S^2$), with adversarial maxima between $0.5$ and $80$, far above the $10^{-6}$ threshold. A recurring feature in these learning runs is that the mean residual can decrease while substantial local defects persist. This behaviour is examined alongside the successful reconstructions on the supplied fibration models. 

\subsection{Plan of paper}
We provide a short review of the  mathematical preliminaries in \S\ref{sec:background},  covering almost complex structures and integrability, and a brief description of the complex manifold benchmarks included in our search. \S\ref{sec:methods} develops the numerical formulation of our programme, and \S\ref{sec:results} presents the outcome of the direct searches. \S\ref{sec:reconstruction} studies the reconstructions of the six-dimensional examples from a geometric input. We conclude the main part of the paper in \S\ref{sec:discussion} with analysis of the opportunities and limitations of the numerical evidence, with the implications for future searches. 

We also include two appendices for supplementary material. Appendix \ref{app:orientation} takes care of implementing an orientation invariant in the searches in \S \ref{sec:results}, and Appendix \ref{app:implementation} records principal settings for these numerical searches.
\subsection*{Acknowledgements}
\textbf{EH} is supported by São Paulo Research Foundation (FAPESP) grant 2024/18994-7. \textbf{TSG} is supported by the 2024
Max Planck-Humboldt Research Award, bestowed on Geordie Williamson by
the Max Planck Society and the Alexander von Humboldt Foundation and
hosted by Catharina Stroppel at the University of Bonn. \textbf{TSRS} is supported by FAPESP grant 2022/09891-4. \textbf{JS} is supported by FAPESP grant 2023/02809-3. \textbf{JS} would also like to thank Robert Bryant and Michael Albanese for their helpful references on complex structures, as well as the Abdus Salam International Centre for Theoretical Physics for their hospitality while carrying out some parts of this paper.    

This research utilised computational resources of the ``\textit{Centro Nacional de Processamento de Alto Desempenho em São Paulo} (CENAPAD-SP)".
This work was supported in part by a grant of access to OpenAI models through the ChatGPT for Academic Researchers program.

\subsection*{Data Availability}
The code and datasets used in this work are available here: \url{https://github.com/edhirst/learning-acs-cs}. 

\subsection*{AI Usage Statement} The methodology and first draft of this work were those of the authors. LLM agents then recommended and implemented revisions which developed ideas and expanded the methodology. These revisions were shepherded and validated by the authors, who take responsibility for the content. The codebase was also developed with assistance from Fable (Anthropic).

\section{Mathematical background}\label{sec:background}
Every complex manifold of complex dimension $n$ is also a smooth real manifold of dimension $2n$, by separating each complex coordinate into its real and imaginary parts. The converse question is associated with the work of Hopf~\cite{Hopf48}: given a smooth even-dimensional manifold, when can its local coordinates be chosen so that the transition maps are holomorphic? This is still open, in general, in real dimension six and above; an overview of this problem, as well as a historical account, is given in~\cite{Agricola:2018hop}.

\subsection{Obstructions to the existence of complex structures}

There are two key obstructions to finding a holomorphic atlas on a given real, compact, even-dimensional  manifold. The first is purely topological: whether the manifold admits an \textit{almost complex structure} (ACS) i.e. a smooth endomorphism $J$ on the tangent bundle that squares to minus the identity. 

For spheres, the classical result is particularly restrictive: among spheres, only $S^2$ and $S^6$ admit almost complex structures~\cite{BS53}. On $S^2$, the usual cross product gives $J_xv=x\times v$. The analogous construction on $S^6\subset\R^7$ uses the octonionic cross product. More generally, Calabi's construction supplies almost complex structures on oriented hypersurfaces in $\R^7$~\cite{Calabi58}. In particular, the topological obstruction to admitting almost complex structures vanishes in all the two and six dimensional cases considered in this paper. 

The second obstruction is a classical theorem of Newlander and Nirenberg~\cite{NN57}: $J$ can be used to define a holomorphic atlas on the manifold if and only if the Nijenhuis tensor $N_J$ vanishes, which is a first order partial differential equation in $J$. For some almost complex structure $J$, $N_J$ is given by the expression
\begin{equation}\label{eq:nijenhuis}
N_J(X,Y):=[JX,JY]-J[JX,Y]-J[X,JY]-[X,Y],
\end{equation}
where $X,Y$ are smooth vector fields and the bracket is the Lie bracket. Although brackets involve derivatives, the expression is tensorial in $X$ and $Y$. If $N_J$ vanishes, we call such a $J$ \textit{integrable}, or simply refer to it as a complex structure.

This separates the problem into two stages. One first seeks a global section of the bundle of linear complex structures, and then solves a first-order nonlinear differential equation within that space. The distinction also provides a natural numerical safeguard: an arbitrary matrix field with a small Nijenhuis-type residual is not necessarily close to a complex structure, because it may fail the algebraic condition. A simple example of this would be the zero field.

In real dimension two, integrability is automatic. Indeed, $N_J$ is alternating and satisfies the relation $N_J(JX,Y)=-JN_J(X,Y)$. For any non-zero $X$, the pair $X,JX$ spans the tangent plane, and since  $N_J(X,JX)=0$, the entire tensor vanishes. Thus $S^2$ and $T^2=S^1\times S^1$ test the construction of almost complex structures, but do not provide a non-trivial test of the integrability solver. In dimensions four and six, the differential equation $N_J=0$ carries additional information.

\subsection{Examples, controls and recent constructions}
The classical positive controls included in this paper are chosen to include examples admitting both K\"ahler and non-K\"ahler metrics. The product $S^2\times S^2$ is the underlying smooth manifold of $\mathbb{CP}^1\times\mathbb{CP}^1$, the product of two projective lines. For $m\ge2$, the quotient
\begin{equation}\label{eq:hopfquotient}
(\C^m\setminus\{0\})/\langle z\mapsto e^{2\pi}z\rangle
\end{equation}
is a Hopf manifold with underlying smooth manifold $S^{2m-1}\times S^1$. This follows by writing $z=e^t x$, where $|x|=1$ and $t\in\R/2\pi\Z$. The cases $m=2,3$ give our four- and six-dimensional Hopf controls. A further example is on the product $S^3\times S^3$, which carries Calabi--Eckmann structures~\cite{CE53}; its Lie group description also relates it to Samelson's construction~\cite{Samelson53}.

 We include three further examples in this paper. The first two are based on recent constructions of complex structures on $S^6$ ~\cite{AC26,Engel26} and $S^4 \times S^2$ \cite{GFL26,LLZ26}; the third is based on an explicit period model for a proposed complex structure on $S^4\times T^2$. We give more details of these constructions as necessary in \S \ref{sec:reconstruction}; the numerical work here neither proves their topological conclusions nor supplies an independent proof of the Hopf problem. The main feature they share is that they all admit a description as a singular torus fibration over $\mathbb{CP}^1$, with exceptional fibres glued in at a finite number of isolated points in $\mathbb{CP}^1$.   
 
Finally, we note here that the metric used to measure a residual should not be confused with a metric constraint on the structure. We use the round metric and its products to define norms, but do not impose $J^\top J=I$ on tangent spaces. This freedom is essential on $S^6$, where LeBrun's theorem excludes complex structures orthogonal for the round metric~\cite{LeBrun87}. Even on the classical examples, an integrable structure need not be orthogonal for the particular background metric used in training.

\section{Neural and polynomial methods}\label{sec:methods}
\subsection{A global ambient representation}
Other similar works \cite{Hirst:2025seh} set up the geometric learning problem by assigning a field to each chart of an atlas and training the fields to agree on the chart overlaps using the atlas transition functions. This method respects the local description of a complex manifold, but introduces a second optimisation problem on top of the ACS and CS goals where local endomorphisms must first define one global tensor. 
This work instead represents the manifolds in terms of their standard embeddings of each constituent sphere, avoiding the need for this extra globally consistent learning goal.

To define the used standard embeddings, write
\begin{equation}\label{eq:embedding}
M=S^{m_1}\times\cdots\times S^{m_k}\subset\prod_{A=1}^k\R^{m_A+1}=\R^N,
\qquad d=\sum_A m_A,\quad N=d+k.
\end{equation}
At $p=(x_1,\ldots,x_k)$, with $|x_A|=1$, the tangent projector is
\begin{equation}\label{eq:projector}
P(p)=\diag\bigl(I-x_1x_1^\top,\ldots,I-x_kx_k^\top\bigr),
\qquad P^2=P=P^\top.
\end{equation}
Any smooth ambient matrix field $M_\theta(p)$ defines a tangent endomorphism by
\begin{equation}\label{eq:compressed}
K_\theta(p)=P(p)M_\theta(p)P(p).
\end{equation}
It vanishes on the normal space and maps the tangent space into itself. Since $P$ restricts to the identity on $TM$, the almost complex condition becomes $K_\theta^2=-P$. As this is a global representation, no overlap loss is required. But it does not remove the topological obstruction on $S^4$, such that every smooth field of the form \cref{eq:compressed} must fail the almost complex identity somewhere.

The compression is deliberately redundant. Distinct ambient matrices can induce the same tangent field, and the identity $|x_A|^2=1$ makes some polynomial features dependent after restriction. These redundancies are manageable at the dimensions considered here and avoid choosing a global tangent frame (which is not always possible). Local orthonormal frames are then used only to evaluate compact residual vectors more efficiently at individual sample points.

\subsection{Intrinsic differentiation and residual norms}
For vector fields in Euclidean coordinates we use $[U,V]=(DV)U-(DU)V$. Tangent fields along $M$ can be extended off the manifold to evaluate this expression; its restriction to $M$ is independent of the extensions. Let $U_i=Pe_i$, where $(e_i)_{i=1}^N$ is the standard ambient basis. For a general endomorphism $K$ on $TM$, satisfying $K=P K P$, define
\begin{equation}\label{eq:ambientN}
\mathcal N_K(U_i,U_j)=[Ke_i,Ke_j]-K[Ke_i,Pe_j]-K[Pe_i,Ke_j]+K^2[Pe_i,Pe_j].
\end{equation}
This is the Nijenhuis torsion of $K$, since $PK=KP=K$. When $K^2=-P$, it agrees with \cref{eq:nijenhuis}. Keeping $K^2$ in the last term makes the formula tensorial even during optimisation, when this algebraic identity may not hold. Therefore the computation uses only $K,P$ and their first ambient derivatives.

\begin{proposition}[Norm from the projected ambient frame]\label{prop:parseval}
For any bilinear map $B:T_pM\times T_pM\to T_pM$ and any orthonormal tangent basis $(\varepsilon_a)_{a=1}^d$,
\begin{equation}\label{eq:parseval}
\sum_{i,j=1}^N\norm{B(Pe_i,Pe_j)}^2
=\sum_{a,b=1}^d\norm{B(\varepsilon_a,\varepsilon_b)}^2.
\end{equation}
\end{proposition}
\begin{proof}
Expand $Pe_i=\sum_a\langle e_i,\varepsilon_a\rangle\varepsilon_a$. After applying bilinearity, each sum over an ambient index contains $\sum_i\langle e_i,\varepsilon_a\rangle\langle e_i,\varepsilon_b\rangle=\delta_{ab}$. The two contractions give the right-hand side.
\end{proof}
\begin{remark} Here, and elsewhere in the paper (unless otherwise indicated), for general tensors we use the pointwise norm associated to the metric induced by the embedding \eqref{eq:embedding}, e.g.  $\|B\|^2 = \sum_{a,b}\|B(\varepsilon_a,\varepsilon_b)\|^2$ for an orthonormal frame $(\varepsilon_a)$; for endomorphisms, this is the Frobenius norm $\|\cdot\|_F$.
\end{remark}
Thus summing \cref{eq:ambientN} over all ordered ambient input pairs gives precisely the tensor norm of the round product metric. A redundant frame does not change the norm. It can, however, increase the number of residual entries seen by a least-squares solver. If $E$ is an $N\times d$ orthonormal tangent frame at the sample point, antisymmetry gives
\begin{equation}\label{eq:intrinsicrows}
\norm{\mathcal N_K}^2=2\sum_{a<b}\norm{E^\top\mathcal N_K(E_a,E_b)}^2,
\qquad
\norm{K^2+P}_F^2=\norm{E^\top K^2E+I_d}_F^2.
\end{equation}
The corresponding residual vector has $d\binom d2+d^2$ entries per point, including the factor $\sqrt2$ in the antisymmetric part. In dimension six this gives $126$ entries. The frame is a device for evaluating the norm at fixed points, not an additional learned geometric field.

\subsection{Sampling and the meaning of numerical validation}\label{sec:validation}
Training expectations use the normalised round product measure. Independent standard Gaussian vectors $g_A\in\R^{m_A+1}$ are normalised to $x_A=g_A/|g_A|$, giving independent uniform samples on each factor. Fresh points are drawn with each new sample batch throughout training and new independent points are drawn for evaluation. For the two main residuals we write
\begin{equation}\label{eq:residuals}
r_K(p)=\norm{K(p)^2+P(p)}_F,
\qquad n_J(p)=\norm{N_J(p)}.
\end{equation}
Then across these residuals means, root-mean-square values, and maxima answer different questions; as stated in the presented results.

Uniform control is especially important for an existence problem. A small mean does not exclude a defect on a set of small volume, and a random sample may miss that set. To locate such regions we use multi-start ascent of the pointwise residuals. The $n$ starting points are those with the largest residual among $32n$ uniform samples. Each ascent step projects the gradient onto the tangent space, moves a fixed distance along the normalised direction and renormalises every sphere factor; the step is kept only if the residual increases. The distance is multiplied by $1.2$ after an accepted step and by $0.5$ after a rejected one, starting at $0.02$ and kept within $[10^{-6},0.3]$. Ascent is used to validate every run. Only the tail-aware conjugation runs of \cref{sec:results} also use it in training: their batches of $192$ points combine $96$ uniform points, the $72$ points with the largest residual among $288$ uniform candidates, and $24$ points from short ascents ($3$ ascents of $8$ starts and $20$ steps), refreshed every $500$ steps. The final ACS checks use $256$ starting points for independent ascent runs, and $150$ ascent steps per run, while the standard Nijenhuis checks use $128$ starts and $75$ steps.

The largest values detected are denoted by $\radv$ and $\Nadv$ respectively. They are empirical lower bounds on the corresponding true suprema. Finding a point with $r_K\ge1$ conclusively invalidates the strict residual test for that field; however failing to find one is not a proof that every point passes. A rigorous global certificate would require additional bounds, for example interval arithmetic or derivative estimates combined with a sufficiently fine covering. This direction is technical, computationally demanding, and specific to each structure, so is not pursued here and left to future developments.

The ambient runs use double precision to represent the numerical values, and evaluate displayed norms as $(\norm{T}^2+10^{-20})^{1/2}$. The regulariser keeps the square root differentiable where a residual vanishes. Consequently, the repeated value $10^{-10}$ is a reporting floor: it means that the residual is below $10^{-10}$. It is not a limit of double precision, which resolves much smaller values (the reconstruction diagnostics of \cref{sec:reconstruction} report residuals of order $10^{-15}$).
Those diagnostics use a separate coordinate implementation, which is another reason to keep the two sets of diagnostics distinct.

\subsection{Learning an almost complex field}\label{sec:acs_learn}
The first objective is simply
\begin{equation}\label{eq:acsloss}
\mathcal L_{\mathrm{ACS}}(\theta)=\E_p\norm{K_\theta^2+P}_F^2.
\end{equation}
We begin with ambient polynomial features, $M_\theta=\sum_f\phi_f W_f$, with independent trainable matrices $W_f$. The feature classes progress from constants and linear coordinates to cross-factor bilinears, homogeneous quadratics and unions of these classes. For three factors, products involving one coordinate from each factor are added at the final stage. These are low-degree polynomial classes, motivated by low-frequency approximation on products of spheres; the sequence is not a fully nested ordering of every harmonic subspace. Each attempt uses $3000$ Adam steps~\cite{Adam}, and the first class passing the preliminary numerical admissibility test is retained. An attempt passes if the largest value of $\norm{K^2+P}$ found by a short worst-point search (\cref{sec:validation}, with $64$ starts of $50$ steps) is below $1$; this threshold is explained in \cref{sec:admissibility}.

There is a reason to start with a small function class. A flexible field can reduce a mean-square loss by making most of the manifold nearly admissible while concentrating the failure in a narrow region. A reduced example on $S^2$ makes the mechanism explicit. Let $J_0v=x\times v$, $f(x)=\langle x,n\rangle$ for a fixed unit vector $n$, and consider $K_{a,b}=(a+bf)J_0$. Since $f$ is uniformly distributed on $[-1,1]$,
\begin{equation}\label{eq:wall}
R(0,b)=\E\norm{K_{0,b}^2+P}^2
=2\left(1-\frac23b^2+\frac15b^4\right).
\end{equation}
At $b^2=5/3$, the Hessian of $R(a,b)$ is positive definite, since its diagonal entries are $16/3$ and its mixed entry vanishes. Thus this two-parameter family has a strict local minimum with opposite complex orientations in the two hemispheres and degeneracy on the equator. This example illustrates a trap in a restricted family, as a local minimum in its function class instead of over the full space of smooth fields.

If plain training of a class fails this test, the class is retried with each of two additional loss terms before the next class is tried. The first adds an annealed Dirichlet term
\begin{equation}\label{eq:dirichlet}
\mathcal E_D(K)=\E\sum_i\norm{P(D_{Pe_i}K)P}_F^2=\E\norm{\nabla K}^2,
\end{equation}
with coefficient decreasing from $3$ to $10^{-3}$. For a tangential endomorphism, $P(D_XK)P$ is its Levi-Civita covariant derivative. The training goal of low losses means low derivatives are then searched for with this term, such that it suppresses an orientation wall (like in the previous simple example) when a low-energy global structure is available, but it is not a universal cure since a non-parallel almost complex structure can itself have significant derivative energy.

The second modification uses an orientation penalty. In a positively oriented orthonormal tangent frame, set
\begin{equation}\label{eq:orientation}
o_K(p)=2^{-d/2}\Pf\bigl(E^\top(K^\top-K)E\bigr),
\end{equation}
where $\Pf$ denotes the Pfaffian of an antisymmetric matrix, the polynomial in its entries satisfying $\Pf(A)^2=\det A$ and $\Pf(B^\top AB)=\det(B)\Pf(A)$; unlike the determinant, its sign changes with the orientation of the frame.
For an exact linear complex structure, $|o_K|\ge1$ and its sign records the induced complex orientation; a proof is given in \cref{app:orientation}. The additional term built from this is $\E\relu(1-so_K)^2$, with $s$ the sign of the mean of $o_K$ over $2000$ points at initialisation, and it was verified to vanish on every exact structure of the selected orientation. It therefore does not specify a target structure, although it does affect the optimisation away from the constraint. Neither additional loss term guarantees that all defects will disappear.

After a successful polynomial attempt, a small-output auxiliary neural network is added as a second learning stage to further reduce the almost complex residual. The matrix field becomes a polynomial plus a multilayer perceptron with GELU activations, four hidden layers of width $256$ in dimension six and three of width $128$ otherwise. Both parts are refined for $10\,000$ Adam steps, followed by up to $200$ iterations of L-BFGS~\cite{LBFGS}, a more accurate but slower second-order optimiser, which are kept only if they lower the held-out loss. The refined field is kept only if it also lowers the worst-point value of the accepted polynomial field (64 starts of 50 steps); otherwise the polynomial field is retained. This happens for $S^2$ (seed 8), $S^1\times S^1$ (seed 7) and $S^3\times S^1$ (seed 7), whose polynomial fields already have worst-point values below $10^{-8}$. These stages finalise the learning of the almost complex structure, ensuring a low residual before the algebraic correction below. More detail about the run settings are collected in \cref{app:implementation}.

\subsection{A conditional correction to the exact algebraic constraint}\label{sec:admissibility}
The threshold used above has a mathematical meaning. If $K|_{T_pM}$ has a real eigenvalue $\lambda$ with unit eigenvector $v$, then $(K^2+P)v=(\lambda^2+1)v$, and hence $r_K(p)\ge1$. Therefore $r_K<1$ excludes real tangent eigenvalues. More directly, the same threshold puts the matrix needed for correction in a domain where its principal inverse square root is well defined.

\begin{proposition}[Algebraic correction]\label{prop:correction}
Let $K=P K P$ be a smooth field on a compact manifold with orthogonal tangent projector $P$. If $\sup_{p\in M}r_K(p)<1$, set
\begin{equation}\label{eq:correction}
S=-K^2+I-P=I-(K^2+P),\qquad J^*=KS^{-1/2},
\end{equation}
where the square root is the principal matrix function. Then $J^*$ is smooth, $P J^*P=J^*$ and $(J^*)^2=-P$.
\end{proposition}
\begin{proof}
The operator norm of $S-I$ is strictly less than one, so the spectrum of $S$ lies in the open disc centred at $1$ with radius $1$. The principal inverse square root exists and varies smoothly~\cite{Higham08}. Since $S$ commutes with $K$ and $P$, so does $S^{-1/2}$. On the normal space $K=0$ and $S=I$; on the tangent space $S=-K^2$. It follows that $(J^*)^2=K^2S^{-1}=-P$, with the stated tangency.
\end{proof}
The proposition is exact, but its hypothesis is a true uniform bound. In the numerical procedure we check $\radv(K)<1$ and then test the corrected field again. This supplies strong practical evidence when a large margin separates the observed residual from one, while retaining the distinction from a theorem about the computed field on all of $M$.

We evaluate the inverse square root by $30$ coupled Newton--Schulz iterations. Starting with $Y_0=S/2$ and $Z_0=I$, set
\begin{equation}\label{eq:newtonschulz}
T_j=\tfrac12(3I-Z_jY_j),\qquad
Y_{j+1}=Y_jT_j,\qquad Z_{j+1}=T_jZ_j.
\end{equation}
Then $KZ_j/\sqrt2$ approaches $J^*$. The strict residual hypothesis implies $\norm{I-S/2}_{\mathrm{op}}<1$, a sufficient convergence condition. Matrix multiplication makes this operation differentiable and compatible with the automatic differentiation used for $N_J$. 
Note that a symmetric eigensolver would be inappropriate here because $-K^2$ need not be symmetric. Also, the correction is algebraic, and does not enforce integrability, so the Nijenhuis tensor must be recomputed after this correction is applied.

\subsection{Conjugation and the neural search for integrability}
Each learnt ACS (after correction) is then used as a starting point for attempted learning of a CS.
The first CS parametrisation used for this learning stage, described in this section, keeps the learned corrected field $J_b$ as a fixed base and varies it by conjugation using a single learnt MLP that outputs both $B_1$ and $B_2$,
\begin{equation}\label{eq:conjugation}
J=gJ_bg^{-1},\qquad g=e^{B_2}e^{B_1},\qquad B_s=P\,\mathrm{MLP}(p)_s\,P,
\end{equation}
Each exponential commutes with $P$ and is invertible, so conjugation preserves $J_b^2=-P$. The loss objective to minimise over the training (via iterative update of the MLP parameters) is $\E\norm{N_J}^2$. General real generators are allowed; restricting to orthogonal conjugation would retain the pointwise singular values of $J_b$ and could unnecessarily constrain the search.

At one point, two exponentials suffice to represent any orientation-preserving real change of frame, where its polar decomposition is a positive definite matrix times an element of $SO(d)$, and both factors are exponentials (parameterised in general by our two learned exponential factors). Globally, however, smoothly choosing logarithms is an additional condition. In fact, a field given by smooth exponential generators is connected to the identity through smooth gauge transformations. The ansatz consequently stays in the homotopy class of $J_b$ and preserves the complex bundle isomorphism type, including its Chern classes. Hence this pointwise expressivity is not equivalent to coverage of every global almost complex structure of a chosen orientation; such that this method has its limitations.

This limitation is concrete on $S^3\times S^1$. Relative to a fixed structure of the same orientation, almost complex structures there carry an integer homotopy invariant in $\pi_3(S^2)=\Z$, the Hopf invariant computed on the slices $S^3\times\{\theta\}$; since $H^2(S^3\times S^1)=0$, Chern classes do not detect it. Relative to the Hopf structure, the learned base of seed 8 has invariant $+1$, so no exponential conjugation of it reaches the Hopf structure (\cref{sec:supervised}). On $S^5\times S^1$, $S^3\times S^3$ and $S^6$ there is a single class for each orientation; on $S^4\times S^2$ and $S^4\times T^2$ classes can differ in $c_1$, and we have not determined the classes of our bases there.

For comparison, additional runs use affine conjugators $g=(I+B_2)(I+B_1)$, with $g^{-1}$ computed explicitly. Although cheaper, this ansatz does require invertibility conditions that exponentials satisfy automatically. Near a point where $g$ becomes singular, $g^{-1}$ and its derivatives become very large, while $J$ stays algebraically exact wherever $g$ is invertible, so the algebraic residual cannot detect the problem. We therefore monitor the smallest singular value $\sigma_{\min}(g)$, on fresh points and by adversarial minimisation, in every conjugation run. Unlike the exponential ansatz, the affine one can also leave the homotopy class of $J_b$.

To reduce differentiation cost, the base field and its derivative are evaluated once per batch. The derivative of the conjugated field is then assembled from
\begin{equation}\label{eq:conjderivative}
DJ=(Dg)J_bg^{-1}+g(DJ_b)g^{-1}-J(Dg)g^{-1}.
\end{equation}
Training again uses Adam followed by held-out-controlled L-BFGS steps for refinement. In the tail-aware runs, batches also include high-residual points (\cref{sec:validation}).

\subsection{Free polynomial fields and second-order optimisation}
The second CS parametrisation returns to a free field, without using MLPs, this field takes a similar form to the ACS case, as
\begin{equation}\label{eq:polyfield}
K_W=P\left(\sum_f\phi_fW_f\right)P,
\qquad
\mathcal L(W)=\E\norm{\mathcal N_{K_W}}^2+10\E\norm{K_W^2+P}^2.
\end{equation}
Here $\phi_f$ contains all ambient monomials of total degree at most $q$. There are $\binom{N+q}{q}$ features and $N^2\binom{N+q}{q}$ coefficients. For example, the degree-two class has $2880$ coefficients when $N=8$. 
For the $S^4 \times T^2$ manifold, which has 3 factors ($S^1, S^1, S^4$), its degree two polynomial is augmented by all possible trilinear monomials (i.e. a linear contribution from each constituent sphere). 
Throughout, degree refers to the degree of the matrix field (inside the brackets of \eqref{eq:polyfield}), and the multiplication by $P$ will then increase the full degree of the resulting ambient expression for $K_W$.

Both $K_W$ and $DK_W$ are linear in $W$. The torsion \cref{eq:ambientN} and the algebraic residual are consequently quadratic in the coefficients. The problem is therefore a structured nonlinear least-squares problem, for which Gauss--Newton and Levenberg--Marquardt methods are natural~\cite{Levenberg44,Marquardt63}. If $R(W)$ is the stacked residual and $Q=\partial R/\partial W$, a damped step solves
\begin{equation}\label{eq:lm}
\left(Q^\top Q+\lambda\diag(Q^\top Q)\right)\Delta W=-Q^\top R,
\end{equation}
to update the learnt weights $W \leftarrow W + \Delta W$.
The damping, controlled by $\lambda$, is reduced after successful steps to prioritise the faster first Gauss-Newton term, then is increased after unsuccessful ones to prioritise the slower but more stable Levenberg–Marquardt term.

Polynomial differentiation is explicit. In particular,
\begin{equation}\label{eq:polyderivative}
\partial_cK=(\partial_cP)MP+P(\partial_cM)P+PM(\partial_cP),
\qquad \partial_cM=\sum_f(\partial_c\phi_f)W_f.
\end{equation}
Quantities independent of $W$ are precomputed on each batch. Also, the residual Jacobian factors through $(K,DK)$, which avoids a separate spatial differentiation for every parameter. These efficiencies are particularly useful when the feature count grows but the ambient dimension remains small.

Each polynomial search starts by projecting its own learned ACS into the chosen class by least squares, matching the classes parameters to the ACS function as best as possible for this second learning stage starting point. This projection may fail the admissibility test of \cref{sec:admissibility}; but admissibility is later checked again on the final candidate. The first class receives $2000$ Adam warm-up steps, and is then optimised in LM rounds, each with a fresh dataset of points and eight iterations. Sets have at least $128$ points and at least twice as many residual entries as coefficients. A held-out loss on fresh points governs stopping and whether the degree is increased. The resampling limits overfitting but cannot rule it out, particularly when a candidate develops narrow defects.

A degree is accepted when its held-out loss is below $10^{-16}$ and the corrected field passes the algebraic checks and the independent condition $\Nadv<10^{-6}$ on the Nijenhuis tensor (64 starts of 60 steps). Within a degree, the LM rounds stop early once the held-out loss falls below $10^{-20}$; this is why accepted runs need not reach the same final loss. Note that the method described in this section changes both the representation and the optimiser relative, compared to neural conjugation of the previous section, so the individual contributions of each cannot be isolated.

There is nevertheless a useful representability check, in confirming the known complex structures are actually inside the polynomial class being used in this sections ansatz. Let $I_{2m}$ be the standard complex structure on $\R^{2m}$, $\xi=I_{2m}x$ on $S^{2m-1}$ and $\eta=J_2y$ on $S^1$. A Hopf structure is obtained by compressing
\begin{equation}\label{eq:hopfmatrix}
M(x,y)=\begin{pmatrix}I_{2m}&\xi\eta^\top\\-\eta\xi^\top&0\end{pmatrix}.
\end{equation}
On $\{x,\xi\}^{\perp}$ it acts as $I_{2m}$, while $J\eta=\xi$ and $J\xi=-\eta$. It is the structure induced by \cref{eq:hopfquotient} in the product coordinates. Similarly, with $\xi_a=I_4x_a$, the Calabi--Eckmann example on $S^3\times S^3$ is represented by
\begin{equation}\label{eq:cematrix}
M(x_1,x_2)=\begin{pmatrix}I_4&-\xi_1\xi_2^\top\\\xi_2\xi_1^\top&I_4\end{pmatrix}.
\end{equation}
These ambient matrices are quadratic; the product structure on $S^2\times S^2$ is linear, so all are reachable within degree two. The checks establish that these classes can contain integrable structures. They are performed separately from training, and no coefficients from these formulas are supplied to the solver.

\section{Learning results}\label{sec:results}

The procedure of \cref{sec:methods} is now applied to the ten spheres and sphere products in this study. The experiments separate the algebraic problem of finding an almost complex structure from the differential problem of making that structure integrable. This distinction is essential, since a small Nijenhuis residual is meaningful for the second problem only when the algebraic constraint has also been controlled. This section is structured such that first the search and approximation of almost complex structures on all ten manifolds is completed, followed by the recovery of classical complex structures where possible, and finally the alternative approaches are presented for the manifolds where the first searches were inconclusive.

Unless stated otherwise, reported means use $6000$ fresh points from the round product measure, and the residual conventions, adversarial diagnostics and $10^{-10}$ reporting floor are those of \cref{sec:validation}. 

\subsection{Almost complex structures}

The almost complex stage passes the numerical admissibility checks (as described in \cref{sec:admissibility}) on both seeds for nine of the ten manifolds. The exception is $S^4$, in agreement with its known topological obstruction. \Cref{tab:acs-paired} records the accepted polynomial class and the residual loss score of the selected field, after optional refinement and before correction (see \cref{sec:acs_learn}). For the accepted fields these residuals range from the $10^{-10}$ reporting floor ($S^1\times S^1$, seed 7) to $0.12$ ($S^4\times T^2$, seed 8). Applying the matrix correction subsequently brings the tested algebraic residual to the $10^{-10}$ reporting floor in every accepted run.

The accepted function classes are usually quite simple. Linear features suffice on $S^2$, $S^6$ and $S^2\times S^2$, while bilinear features suffice on $S^1\times S^1$, $S^3\times S^3$, $S^5\times S^1$ and $S^4\times S^2$. On $S^3\times S^1$ the two seeds select different classes, linear and homogeneous quadratic. The three-factor case $S^4\times T^2$ required these trilinear features. Note these aren't necessarily the minimal-degree polynomial representation, just the outcomes of this search process in these set function classes.

\begin{figure}[tb]
\centering
\includegraphics[width=0.98\textwidth]{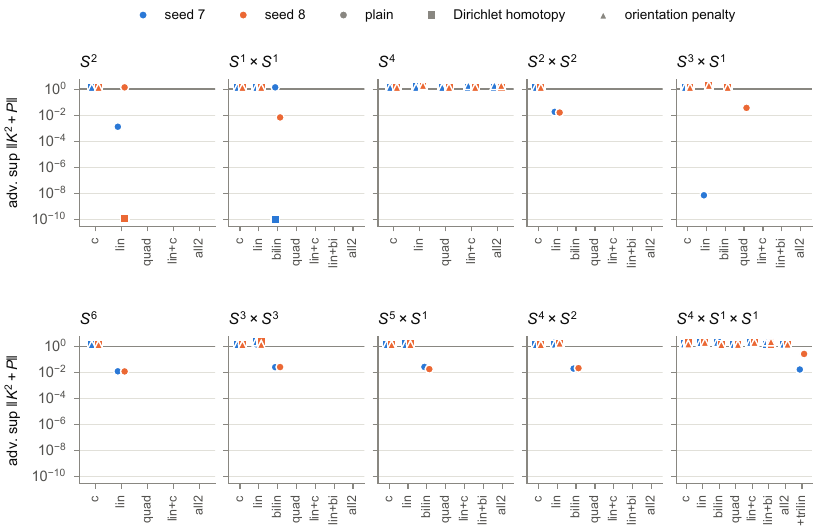}
\caption{The almost complex feature ladder for each considered manifold for both seed runs (the two colours). Circles, squares and triangles denote plain training, the Dirichlet addition and the orientation-penalty addition respectively. The horizontal line is the sufficient admissibility threshold $\norm{K^2+P}<1$; the plotted quantities are the preliminary ascent maxima (64 starts of 50 steps) of each attempt's polynomial field, before refinement. Each point is one attempt; the function class is made progressively more complex until the threshold is met. The selected final fields, with refinement retained according to \cref{sec:acs_learn}, are presented in \Cref{tab:acs-paired}. The function classes, in the order tried, are: c, constant; lin, the ambient coordinates; bilin, products of two coordinates from different sphere factors; quad, all quadratic monomials; lin+c and lin+bi, the unions of lin with c and with bilin; all2, all monomials of degree at most two; +trilin, all2 together with the products of one coordinate from each of the three factors. In \Cref{tab:acs-paired}, c, lin, bilin, quad and +trilin are denoted C, L, B, Q and T. Classes that need several factors appear only on products.}
\label{fig:results-acs-ladder}
\end{figure}

\begin{table}[tb]
\centering\small
\caption{Residual loss scores for the almost complex searches, shown for the two seeded runs with seeds $\{7,8\}$. The selected class is labelled by the required polynomial degree needed to reach the stated approximation level. Using notation L for linear features, B for products of coordinates from different sphere factors, Q homogeneous quadratic features, and T all features of degree at most two together with three-factor products. C is the constant class retained after rejection of all these classes on $S^4$. Residuals which are numerically below the admissibility threshold are marked Pass. For $S^2$ (seed 8), $S^1\times S^1$ (seed 7) and $S^3\times S^1$ (seed 7) the network refinement did not lower the worst-point residual and was discarded (\cref{sec:acs_learn}); their entries are those of the polynomial field, and $10^{-10}$ is the reporting floor.}
\label{tab:acs-paired}
\begin{tabular}{@{}lccccc@{}}
\toprule
& \multicolumn{2}{c}{Selected class} & \multicolumn{2}{c}{$\radv$} & \\
\cmidrule(lr){2-3}\cmidrule(lr){4-5}
Manifold & seed 7 & seed 8 & seed 7 & seed 8 & Outcome\\
\midrule
$S^2$ & L & L & $4.2\times10^{-5}$ & $1.3\times10^{-10}$ & Pass\\
$S^1\times S^1$ & B & B & $1.0\times10^{-10}$ & $6.8\times10^{-5}$ & Pass\\
$S^4$ & C & C & $2.00$ & $2.00$ & Reject\\
$S^2\times S^2$ & L & L & $9.8\times10^{-5}$ & $8.9\times10^{-5}$ & Pass\\
$S^3\times S^1$ & L & Q & $7.1\times10^{-9}$ & $1.6\times10^{-3}$ & Pass\\
$S^6$ & L & L & $1.3\times10^{-4}$ & $1.5\times10^{-4}$ & Pass\\
$S^3\times S^3$ & B & B & $2.2\times10^{-4}$ & $1.9\times10^{-4}$ & Pass\\
$S^5\times S^1$ & B & B & $9.7\times10^{-5}$ & $1.9\times10^{-4}$ & Pass\\
$S^4\times S^2$ & B & B & $3.0\times10^{-4}$ & $1.6\times10^{-4}$ & Pass\\
$S^4\times T^2$ & T & T & $8.0\times10^{-4}$ & $0.12$ & Pass\\
\bottomrule
\end{tabular}
\end{table}

The reruns with the additional loss terms described in \cref{sec:acs_learn} also reveal sensitivity to optimisation. Adding the Dirichlet homotopy leads to success on $S^2$ for seed 8 and on the torus for seed 7 after plain training has stalled. Across all the manifolds, no run using the orientation penalty was the first successful attempt for each campaign. These results hence demonstrate the usefulness of this ladder of increasing function class complexity, and highlight that plain training is often sufficient, while these additions still do not break the non-existence on $S^4$. The performances of the numerical runs at each function class complexity stage through the ladder are shown in \cref{fig:results-acs-ladder}.

On $S^4$, the mean algebraic residual after refinement is only $0.0106$ for seed 7 and $0.0125$ for seed 8, while the adversarial value is $2.00$ for both, the value $\norm{P}_F$ taken where $K$ vanishes. This contrast demonstrates why a mean alone is inadequate, as the field can satisfy the algebraic equation rather well on average over most of the manifold, but retains a substantial local defect that can't be removed. Both seed runs for this manifold are rejected and their complex-structure stages are therefore not run, since there is no ACS to run from. The topological obstruction to existence explains why a genuine solution cannot be obtained here, and these numerical results provide satisfying corroboration.

The corrected almost complex fields need not be close to integrable. For example, the mean Nijenhuis norm on $S^6$ is approximately $19.60$, close to the value $\sqrt{384}\approx19.60$ of the octonionic structure (\cref{app:implementation}), and it is this norm that will be minimised in the next learning stage, where these learnt almost complex structures are used as the respective starting points.

\subsection{Recovery of classical complex structures}

The polynomial search with Levenberg--Marquardt is successful on both seeds for all classical examples $\{S^2, S^1 \times S^1, S^2 \times S^2, S^3 \times S^1, S^3 \times S^3, S^5 \times S^1\}$, as shown by the low losses in \Cref{tab:classical-learning}. 
Note the already learnt ACS's on the 2d surfaces $\{S^2, S^1 \times S^1\}$ are automatically integrable, but are still included as comparative controls.  
The loss shown is the held-out value of $\E\norm{\mathcal N_K}^{2}+10\E\norm{K^2+P}^{2}$ before correction; the Nijenhuis column evaluates the true integrability condition on the corrected field.

\begin{table}[tb]
\centering\small
\caption{Polynomial complex-structure searches, run for 2 seed runs through the ladder of increasing function class complexity, with the degree of the polynomial in the ambient $\mathbb{R}^N$ coordinates listed. Time is the recorded elapsed time of the CS stage, excluding its preceding ACS stage. Every row passes the numerical acceptance tests, with $\Nadv$ at or near the $10^{-10}$ reporting floor (at most $6.5\times10^{-10}$).}
\label{tab:classical-learning}
\begin{tabular}{@{}lccccc@{}}
\toprule
Manifold & Seed & Degree & Held-out loss & $\Nadv$ & Minutes\\
\midrule
$S^2$ & 7 & 2 & $7.3\times10^{-21}$ & $1.0\times10^{-10}$ & $0.29$\\
& 8 & 2 & $8.2\times10^{-21}$ & $1.0\times10^{-10}$ & $0.29$\\
$S^1\times S^1$ & 7 & 2 & $8.1\times10^{-22}$ & $1.0\times10^{-10}$ & $0.35$\\
& 8 & 2 & $8.3\times10^{-18}$ & $1.0\times10^{-10}$ & $0.35$\\
$S^2\times S^2$ & 7 & 2 & $7.5\times10^{-30}$ & $1.0\times10^{-10}$ & $1.02$\\
& 8 & 2 & $7.4\times10^{-30}$ & $1.0\times10^{-10}$ & $0.95$\\
$S^3\times S^1$ & 7 & 3 & $5.6\times10^{-21}$ & $1.5\times10^{-10}$ & $2.11$\\
& 8 & 2 & $5.8\times10^{-20}$ & $6.5\times10^{-10}$ & $2.19$\\
$S^3\times S^3$ & 7 & 2 & $3.8\times10^{-21}$ & $1.1\times10^{-10}$ & $5.51$\\
& 8 & 2 & $9.6\times10^{-21}$ & $1.1\times10^{-10}$ & $5.57$\\
$S^5\times S^1$ & 7 & 2 & $3.0\times10^{-20}$ & $2.5\times10^{-10}$ & $10.22$\\
& 8 & 2 & $1.4\times10^{-20}$ & $1.7\times10^{-10}$ & $8.95$\\
\bottomrule
\end{tabular}
\end{table}

\begin{figure}[tb]
\centering
\includegraphics[width=0.98\textwidth]{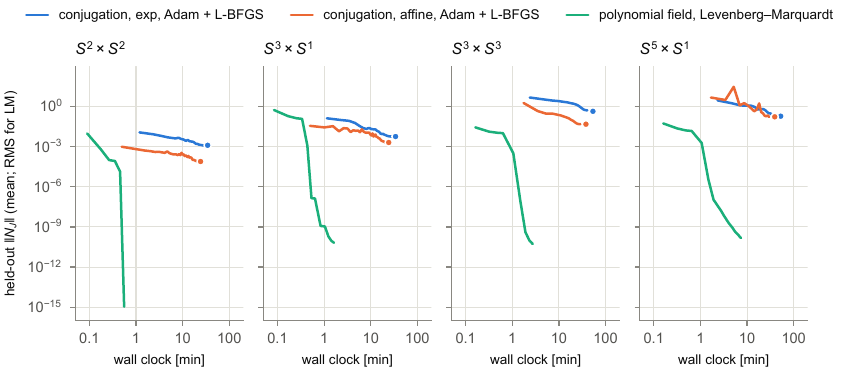}
\caption{Classical examples, seed 7: conjugation-based neural searches and polynomial fields with Levenberg--Marquardt. The neural curves show held-out mean Nijenhuis norms; the polynomial curves show the root mean square of the free-field Nijenhuis norm. These are different summaries, with the root mean square bounding the mean from above. Final corrected-field results are given in \cref{tab:classical-learning}. On $S^5\times S^1$ the affine conjugator is numerically singular (\cref{sec:supervised}), so its curve is not a valid result.}
\label{fig:results-cs-comparison}
\end{figure}

Eight of the 12 runs reach a held-out loss below $10^{-20}$, three more stop between $10^{-20}$ and $6\times10^{-20}$, and the last ($S^1\times S^1$, seed 8) stops at $8.3\times10^{-18}$; all pass the acceptance test of \cref{sec:methods}. Eleven runs stop at degree two. Only $S^3\times S^1$ with seed 7 steps up to degree three, while the same manifold succeeds at degree two with seed 8. The training histories of the neural and polynomial approaches are compared in \cref{fig:results-cs-comparison}.

By comparison, the conjugation methods involving the neural architectures are far inferior. With exponential conjugators, the final mean Nijenhuis norms on $S^2\times S^2$, $S^3\times S^1$, $S^3\times S^3$ and $S^5\times S^1$ are respectively $0.00127$, $0.00574$, $0.441$ and $0.191$. Affine conjugation improves the first three to $7.8\times10^{-5}$, $0.00205$ and $0.048$, but still does not reach the accuracy of the polynomial solver. On $S^5\times S^1$ it reaches $0.171$, but only through a numerically singular conjugator (adversarial $\sigma_{\min}(g)=2.9\times10^{-8}$), so this is not a valid improvement. Additionally, these runs take approximately 24--54 minutes, compared with 1--10 minutes for the corresponding polynomial runs. 

The learning process specifies no specific complex structure to model, such that any solution to the integrability equation is sufficient.
These learnt complex structures can hence be shown to not match the classical examples. The classical product structures, and the Hopf structures with the normalisation of \cref{eq:hopfquotient}, are orthogonal for the round metric: $J^\top J=I$ on tangent spaces, so all singular values of $J$ equal $1$.
For the learnt four- and six-dimensional structures, the mean of $\norm{J^\top J-I}$ ranges from $0.054$ to $4.61$, with singular values between $0.41$ and $2.42$, far from the classical behaviours and hence not the explicit round-orthogonal representatives (but they may still be biholomorphic to them). On $S^3\times S^1$, the complex structure learned from seed 8 is not even homotopic to the Hopf structure through almost complex structures: its Hopf invariant relative to it is $+1$.
Additionally, for some structures a near constant value of this norm may suggest some additional symmetry.

Importantly, the resampling between Levenberg--Marquardt optimisation rounds prevents overfitting of the complex structure to a specific subset of points and focusing the error in a small region outside of the sample.
However, it is still possible that a sufficiently small defect area may not be sampled at all over the training period, which is an unavoidable possibility in numerics, but the independent adversarial checks make it less likely to go unnoticed.

\subsection{A supervised benchmark on the Hopf manifolds}\label{sec:supervised}

The architectures above are trained to minimise the integrability equations alone, with a fixed computational budget. To compare the training signals, the same networks are also trained in a supervised manner using an explicit complex structure. Choosing $S^3\times S^1$ and $S^5\times S^1$ as representatives in complex dimensions 2 and 3, both these manifolds exhibit a complex structure in closed form, namely a Hopf structure.\footnote{We used the Hopf structure of $(\C^m\setminus\{0\})/\langle z\mapsto e\,z\rangle$. It differs from \cref{eq:hopfquotient} only in the scaling factor, but it is far from orthogonal for the round metric, which makes it harder to learn.
} 
The same network architectures (same ansatz form and same size) with exponential conjugation setup where they take input of a point on the manifold and output a prediction of the structure $J_{pred}$, were trained to reproduce this Hopf structure $J_{Hopf}$, in two ways: from its values alone via minimising the supervised loss $|J_{pred}-J_{Hopf}|$, and from its values together with its first derivatives via also minimising $|\partial_iJ_{pred}-\partial_iJ_{Hopf}|$.
\Cref{tab:supervised} provides results for all four training run styles, including for these two supervised runs and the previous two unsupervised runs, comparing through the mean Nijenhuis norm which vanishes for an exact complex structure.
The supervised fits conjugate a learned base field as in \cref{eq:conjugation}, so the base must lie in the homotopy class of $J_{Hopf}$. Bases with the opposite orientation are first reflected ($\theta\mapsto-\theta$) or negated, and on $S^3\times S^1$ we use the base of seed 7 and an additional base from seed 10, whose Hopf invariants relative to $J_{Hopf}$ vanish. From the base of seed 8, whose invariant is $+1$, the fit with derivatives stalls at a mean Nijenhuis norm of $1.58$ (adversarial maximum $1.8\times10^{3}$). The budgets are not matched: the supervised fits use $40\,000$ ($S^3\times S^1$) or $80\,000$ ($S^5\times S^1$) Adam steps followed by L-BFGS, taking from 14 minutes to 19 hours, while the unsupervised runs use $4000$--$10\,000$ Adam steps and 24--53 minutes.

\begin{table}[tb]
\centering\small
\caption{Supervised and unsupervised training on the Hopf manifolds, compared through mean and adversarial Nijenhuis norms; the computational budgets are not matched. Supervised networks are trained to reproduce a known Hopf structure, from its values alone or from its values and first derivatives. Unsupervised networks are those of \cref{fig:results-cs-comparison}, trained on the equations alone. Ranges are over independent runs: bases 7 and 10 on $S^3\times S^1$, and bases 7 and 8 with two network seeds each on $S^5\times S^1$. Supervised means use 4000 fresh points, unsupervised means 6000; $\Nadv$ is the adversarial maximum. $^\dagger$Not a valid result: the conjugator is numerically singular (adversarial $\sigma_{\min}(g)=2.9\times10^{-8}$).}
\label{tab:supervised}
\begin{tabular}{@{}llcc@{}}
\toprule
Manifold & Training & Mean $\norm{N_J}$ & $\Nadv$\\
\midrule
$S^3\times S^1$ & supervised, values only & $0.035$ & $0.10$\\
 & supervised, values and derivatives & $0.015$--$0.016$ & $0.054$\\
 & unsupervised, exponential conjugation & $0.0057$ & $0.016$\\
 & unsupervised, affine conjugation & $0.0021$ & $0.0060$\\
\addlinespace
$S^5\times S^1$ & supervised, values only & $0.093$ & $0.38$\\
 & supervised, values and derivatives & $0.025$--$0.028$ & $0.076$\\
 & unsupervised, exponential conjugation & $0.19$ & $109$\\
 & unsupervised, affine conjugation$^\dagger$ & $0.17$ & $5.2\times10^{15}$\\
\bottomrule
\end{tabular}
\end{table}

Supervision on values alone sets the benchmark. These networks reproduce the known structure closely on average, with relative value errors of $2.3\times10^{-4}$ on $S^3\times S^1$ and $5.1\times10^{-4}$ on $S^5\times S^1$, yet their Nijenhuis norm stays comparatively large. Integrability is a condition on derivatives, and matching values does not control them. On $S^3\times S^1$, the unsupervised networks with exponential conjugation never see the solution, yet they beat this benchmark by a factor of about six: there, the equations themselves are a more effective training signal for integrability than knowledge of the solution's values. On $S^5\times S^1$ this is reversed, and the exponential run is about twice the values-only benchmark.

Supervision with derivatives is a stronger reference. Given both the values and the first derivatives of the solution, the networks reach mean Nijenhuis norms of $0.015$--$0.016$ on $S^3\times S^1$ and $0.025$--$0.028$ on $S^5\times S^1$, about two and three to four times below the values-only fits. On $S^3\times S^1$ both unsupervised networks go beyond even this reference: exponential conjugation by a factor of about three, and affine conjugation, at $0.0021$, by a factor of seven to eight, the best result in \cref{tab:supervised}. Trained on the equations alone, they attain smaller Nijenhuis residuals than networks given the explicit solution together with its derivatives. On $S^5\times S^1$ they do not: both supervised fits beat the exponential run, by factors of about two (values only) and seven to eight (with derivatives), and the affine value is not valid, because its conjugator is numerically singular.

Two remarks qualify the comparison. On $S^5\times S^1$, the exponential run keeps a small region where the Nijenhuis norm is large (adversarial maximum $109$, with the worst 1\% of points carrying 98\% of $\E\norm{N_J}^2$), whereas the supervised fits show no such region (adversarial maxima at most $0.38$). And the polynomial fields of \cref{tab:classical-learning} remain far more accurate than every neural network in this comparison, supervised or not.

\subsection{The unresolved numerical searches}

The same polynomial procedure does not produce an admissible near-integrable field on $S^6$, $S^4\times S^2$ or $S^4\times T^2$. On $S^6$, the returned degree-three field has loss $5.86$ and $\radv=2.48$. On $S^4\times S^2$, the retained degree-two field has loss $6.35$ and $\radv=2.00$; degree three was also tried, with the worse held-out loss $8.58$. On $S^4\times T^2$, the enriched three-factor class has loss $10.65$ and $\radv=3.05$. All three fail the sufficient admissibility test, and hence have not numerically approximated complex structures.

Conjugation preserves the algebraic equation wherever the conjugator is invertible, and exponential conjugation builds in that invertibility. It therefore permits a more direct examination of the remaining differential error. Selected results appear in \cref{tab:hard-learning}. The mean and the adversarial statistic often rank methods differently. On $S^6$, tail-aware sampling lowers the adversarial value from $6.75$ to $3.83$ while increasing the mean from $0.375$ to $1.03$. The same trade-off is stronger on $S^4\times S^2$, where the adversarial value falls from $79.4$ to $10.5$. A polynomial generator does not improve on this ($20.9$) and remains far from the successful classical runs.

\begin{table}[tb]
\centering\small
\caption{Selected six-dimensional searches, seed 7. Means are final fresh-sample estimates. The first row of each block is the corrected ACS starting field. ``Tail'' denotes tail-aware sampling; the continued $S^4\times T^2$ run adds 4000 steps to the initial 1500. The affine row is diagnostic: its conjugator is numerically singular (adversarial $\sigma_{\min}(g)=3.8\times10^{-9}$).}
\label{tab:hard-learning}
\begin{tabular}{@{}llcc@{}}
\toprule
Manifold & Route & Mean $\norm{N_J}$ & $\Nadv$\\
\midrule
$S^6$ & ACS starting field & $19.60$ & $19.60$\\
& Exponential MLP & $0.375$ & $6.75$\\
& Exponential MLP, tail & $1.03$ & $3.83$\\
& Polynomial generator, degree 3 & $1.33$ & $9.78$\\
\addlinespace
$S^4\times S^2$ & ACS starting field & $15.18$ & $16.14$\\
& Exponential MLP & $0.855$ & $79.36$\\
& Affine MLP & $0.168$ & $2.58\times10^{16}$\\
& Exponential MLP, tail & $2.23$ & $10.50$\\
& Polynomial generator, degree 2 & $3.61$ & $20.86$\\
\addlinespace
$S^4\times T^2$ & ACS starting field & $2.98$ & $3.93$\\
& Exponential MLP, tail & $0.138$ & $0.703$\\
& The same, continued & $0.0924$ & $0.479$\\
\bottomrule
\end{tabular}
\end{table}

The affine $S^4\times S^2$ run is especially instructive. Its mean is $0.168$ and its maximum over 6000 points is $18.4$, whereas ascent finds $2.6\times10^{16}$. The conjugator is numerically singular: $\sigma_{\min}(g)$ is $0.21$ on 6000 points but $3.8\times10^{-9}$ under adversarial search. Since an affine parametrisation does not supply the global invertibility guarantee of an exponential, this run cannot be treated as a reliable global improvement merely because its mean is small.

Independent fresh-sample analysis confirms concentration of the differential defect. In the exponential mean-square runs, the worst 5\% of points account for approximately 49\% of the squared residual on $S^6$ and 61\% on $S^4\times S^2$. On $S^5\times S^1$, despite the existence of a classical complex structure, the worst 1\% account for 98\% on a sample of 6000 points. Thus concentration is not specific to the difficult topological cases. Tail-aware training reduces the worst-5\% shares to 17\% and 15\% on the first two manifolds; it distributes the remaining error more evenly without eliminating it.

\begin{figure}[tb]
\centering
\includegraphics[width=0.72\textwidth]{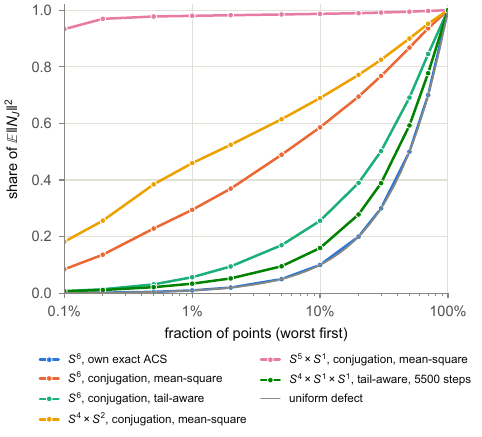}
\caption{Defect concentration measured on fresh samples: the fraction of the total squared Nijenhuis residual carried by the points with the largest residuals. The undotted line represents a uniform contribution from all points, it acts as a baseline along with the exact ACS on $S^6$ (which overlaps with it).}
\label{fig:results-concentration}
\end{figure}

The statistical mechanism is straightforward. If most of the manifold has a modest residual and a small region has a much larger one, a fresh batch may contain few points from that region or none at all. Decreasing its volume can then improve the mean while leaving a serious local failure, as demonstrated in \Cref{fig:results-concentration}. Ascent and tail-aware sampling address this sampling imbalance by deliberately revisiting large-residual regions and this process identified the local defects in these cases. Consequently, a smaller adversarial value is a useful comparative diagnostic, but it does not establish a global upper bound or rule out an undiscovered concentration elsewhere since these finitely many adversarial searches may not have had time to locate a local defect in the learnt structure.

The recent complex structures produced on these manifolds have very specific algebraic constructions, and perhaps not surprisingly this general numerical search struggled to find them within the computational budget.
More focused numerical representation is described in \cref{sec:reconstruction}.

\section{Reconstruction from geometric input}\label{sec:reconstruction}

The searches in \cref{sec:results} ask whether integrability can be recovered from a completely general ansatz on an embedded manifold. A different experiment begins with a geometric description in which much of integrability is already built in. We consider both approaches because their contrast identifies what the optimisation is actually being asked to learn. On a torus fibration, holomorphic period data provide local complex coordinates directly; the remaining numerical problem concerns non-degeneracy and compatibility between those local descriptions. A successful reconstruction therefore measures the effectiveness of a numerical representation informed by geometry. 

For $S^6$ and $S^4 \times S^2$, the input is the torus-fibration construction described in the recent manuscripts \cite{AC26,Engel26,GFL26}; a different construction of $S^4 \times S^2$ appears in \cite{LLZ26}. We take the monodromy, the form of the periods, the compactification data and the stated topological identifications from those works. Our calculation solves the resulting functional equations and checks their numerical consistency. At the time of completing this work there was no published construction on $S^4\times T^2$; therefore one was designed and this input is instead the explicit period model of \cref{sec:s4t2-periods}. The numerical tests in this section check its displayed formulas, but do not by themselves show that the compactified fibration exists or is diffeomorphic to $S^4\times T^2$. This is subject of ongoing work.

\subsection{Period coframes and local integrability}

Let $Z$ be a complex coordinate on the base and let $x\in\R^4/\Lambda$ denote real coordinates on a torus fibre. A complex $2\times4$ matrix $\Pi(Z)$ defines an almost complex structure on the regular part of the family when the forms
\begin{equation}\label{eq:period-coframe-new}
 dZ,\qquad \Pi(Z)\,dx
\end{equation}
together with their conjugates form a complex coframe. Equivalently, the columns of $\Pi$ must be linearly independent over $\R$ in $\C^2$. Writing $\Theta$ for the matrix of the three forms in real coordinates, the induced real endomorphism is determined by $\Theta J=i\Theta$, or explicitly by
\begin{equation}\label{eq:coframe-J-new}
 J=\begin{pmatrix}\Theta\\\overline\Theta\end{pmatrix}^{-1}
   \begin{pmatrix}i\Theta\\-i\overline\Theta\end{pmatrix}.
\end{equation}
The identity $J^2=-I$ follows algebraically wherever the coframe is non-degenerate.

For the first two reconstruction experiments, the normalisation is
\begin{equation}\label{eq:normalised-period-new}
 \Pi=\begin{pmatrix}6\mu&\tau&1&0\\ \beta&\mu&0&1\end{pmatrix},
 \qquad \Im\tau>0.
\end{equation}
In this form non-degeneracy is equivalent to
\begin{equation}\label{eq:period-determinant-new}
 D:=\Im\beta-\frac{6(\Im\mu)^2}{\Im\tau}\ne0.
\end{equation}
The intended component of the period domain has $D<0$. The assumption $\Im\tau>0$ is part of the construction; without it, one should use the undivided determinant condition $\Im\tau\Im\beta-6(\Im\mu)^2\ne0$.

In the normalisation \eqref{eq:normalised-period-new}, integrability is equivalent to holomorphy of $\Pi$. Indeed, closure of the ideal generated by \eqref{eq:period-coframe-new} under exterior differentiation requires $\partial_{\overline Z}\Pi=L\Pi$ for a complex $2\times2$ matrix $L$. The last two columns of $\Pi$ are constant and form the identity, so these columns force $L=0$. Conversely, when $\Pi$ is holomorphic, differentiating its rows gives only terms containing $dZ$, and the ideal is closed. For a general, unnormalised period coframe, holomorphy remains sufficient, which is the implication used below for $S^4\times T^2$.

This elementary observation is important for interpreting the experiment. Once holomorphy is imposed by the representation, a small computed Nijenhuis tensor checks the implementation of \eqref{eq:coframe-J-new} and its derivatives. It is the compatibility residual, rather than the local Nijenhuis residual, that measures the success of the period solve.

\subsection{A spectral reconstruction for \texorpdfstring{$S^6$ and $S^4\times S^2$}{S6 and S4 x S2}}

The construction uses the triangle group $\Delta(3,4,\infty)$. In the cusp coordinate employed by the calculation, the parabolic generator is $g_0Z=Z-1$, while the two elliptic generators act by the matrices
\begin{equation}\label{eq:triangle-new}
 g_1=\begin{pmatrix}1&1-\sqrt2\\1+\sqrt2&0\end{pmatrix},
 \qquad
 g_2=\begin{pmatrix}0&1-\sqrt2\\1+\sqrt2&\sqrt2\end{pmatrix}.
\end{equation}
The periods satisfy the six equations
\begin{align}
 \tau(g_1Z)&=\frac{\tau-1}{\tau},&
 \tau(g_2Z)&=-\frac1\tau,\label{eq:period-tau-new}\\
 \mu(g_1Z)&=\frac{1-\mu}{\tau},&
 \mu(g_2Z)&=1+\frac\mu\tau,\label{eq:period-mu-new}\\
 \beta(g_1Z)&=\beta+2-\frac{6(1-\mu)^2}{\tau},&
 \beta(g_2Z)&=\beta-3-\frac{6\mu^2}{\tau}.
 \label{eq:period-beta-new}
\end{align}
All functions on the right-hand sides are evaluated at $Z$. These equations and their cusp asymptotics are geometric input to the reconstruction, rather than equations inferred by the solver.

With $q=e^{2\pi iZ}$, the numerical representation is
\begin{equation}\label{eq:qseries-new}
 \tau=Z+\sum_{n=0}^{K}a_n(q/q_0)^n,\qquad
 \mu=\sum_{n=0}^{K}b_n(q/q_0)^n,\qquad
 \beta=-Z+\sum_{n=0}^{K}c_n(q/q_0)^n.
\end{equation}
Here $q_0\simeq0.156$ scales the basis to the range encountered on the working domain. Every truncation is holomorphic and satisfies the parabolic transformation law exactly. The nonlinear solve therefore concerns only \eqref{eq:period-tau-new}; after $\tau$ is fixed, the equations for $\mu$ and $\beta$ are linear least-squares problems. The additive constant in $\beta$ is free in the compatibility equations and is fixed at $c_0=-3i$. The recorded tests place the resulting periods in the component $D<0$.

\begin{figure}[t]
 \centering
 \includegraphics[width=\textwidth]{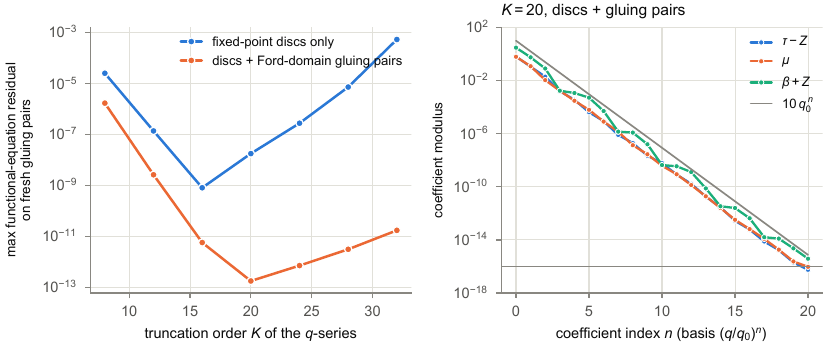}
 \caption{Resolution study for the period solve. Left: largest functional-equation residual on held-out gluing pairs, comparing fixed-point collocation with the addition of gluing pairs. Right: coefficient magnitudes in the scaled basis. These are numerical convergence diagnostics; holomorphy is imposed by the representation, whereas global compatibility is satisfied only to the measured accuracy.}
 \label{fig:recon-convergence-new}
\end{figure}

The collocation set combines small circles about the two elliptic fixed points with pairs of points related by the generators near the domain identifications. This second set is essential. Collocation near the fixed points alone constrains the local functional equations accurately, but leaves combinations of higher coefficients poorly determined near the edges. In the resolution study of \cref{fig:recon-convergence-new}, adding gluing pairs (left panel) reduces the largest held-out functional-equation residual to approximately $1.8\times10^{-13}$ at $K=20$. Larger truncations do not improve on this: the residual rises slowly, to $1.7\times10^{-11}$ at $K=32$, consistent with the coefficients beyond $n\approx20$ having reached the rounding level (right panel). With fixed-point collocation alone, the best value is $8.2\times10^{-10}$ at $K=16$. This convergence is therefore evidence for a useful numerical resolution, rather than a rigorous bound on the omitted series tail.

The lattice action is also checked independently of the period fit. With the monodromy convention of \cite{AC26}, the real fibre coordinates transform by $A_j=(T_j^{-1})^{\mathsf T}$. The corresponding condition on the coframe is $\Pi(g_jZ)A_j=L_j(Z)\Pi(Z)$ for an invertible complex matrix $L_j$. Alternative transpose and inverse conventions give discrepancies of order one or larger. For $S^4 \times S^2$, the prescribed index-two superlattice changes the period matrix to $\Pi B$, where $B=\operatorname{diag}(1,1,1,1/2)$; the transformed matrices $B^{-1}A_jB$ are integral. Thus the same learned functions support both numerical reconstructions, with different lattice data.

\subsection{What the validation establishes}

The numerical checks separate algebraic exactness, local integrability, deck compatibility and blending across an identified edge. For blending, two transported descriptions are combined with a smooth weight and then corrected algebraically to an almost complex structure as in \cref{sec:methods}. Its Nijenhuis tensor detects errors introduced by the transition. This is a useful overlap diagnostic that is indicative of good numerical approximation, but does not lead to a rigorous global atlas (no theorem correcting the approximate field to a truly integrable one is provided). \Cref{tab:reconstruction-new} collects these diagnostics for both models. The algebraic and regular-part Nijenhuis residuals are at the level of rounding error (at most $1.5\times10^{-14}$), while the functional-equation, deck and blended-overlap errors are larger but remain below $6\times10^{-12}$.

\begin{table}[t]
 \centering\small
 \caption{Selected reconstruction diagnostics from the saved runs. Entries are observed maxima, except for the final row, which is a detected minimum. Norms use the coordinates of each model. Nijenhuis and transition errors are distinct quantities. The $S^4\times T^2$ rows use a common filling annulus ($\varepsilon=0.017$).}
 \label{tab:reconstruction-new}
 \begin{tabular}{@{}lcc@{}}
 \toprule
 Quantity & $S^6$ model & $S^4 \times S^2$ model\\
 \midrule
Algebraic residual $\norm{J^2+I}$ & $2.6\times10^{-15}$ & $4.1\times10^{-15}$\\
Regular-part $\norm{N_J}$, adversarial & $9.9\times10^{-15}$ & $1.5\times10^{-14}$\\
Functional-equation residual & $1.8\times10^{-13}$ & $1.8\times10^{-13}$\\
Largest recorded deck discrepancy & $4.5\times10^{-12}$ & $5.6\times10^{-12}$\\
Blended-overlap $\norm{N_J}$ & $4.2\times10^{-12}$ & $5.8\times10^{-12}$\\
\midrule
Quantity & Supplied $S^4\times T^2$ periods & Learned periods\\
\midrule
Regular-part $\norm{N_J}$ & $3.5\times10^{-14}$ & $1.1\times10^{-13}$\\
Largest chart-transition discrepancy & $7.1\times10^{-12}$ & $6.6\times10^{-12}$\\
Normalised determinant, smallest found & $0.626$ & $0.941$\\
\bottomrule
 \end{tabular}
\end{table}

The compactification requires further analytic input. In the cusp coordinate, the matrix
\begin{equation}
 C(q)=\begin{pmatrix}6\mu&\tau-Z\\\beta+Z&\mu\end{pmatrix}
\end{equation}
is holomorphic at $q=0$ for the truncated representation. The supplied checks evaluate its periodicity and the quantities entering the toric filling construction \cite{AC26,Mumford72}. The central fibre itself is governed by that analytic construction. Likewise, logarithmic transforms require their specified holomorphic translations and twisting data; invariance under constant translations in $x$ does not by itself justify a translation depending on the base. We therefore retain the distinction between regular-coordinate tests and the external results used to extend and identify the compactified model.

All norms in \cref{tab:reconstruction-new} are coordinate norms on the fibration models. No explicit diffeomorphism transporting these structures to the round embedded manifolds is computed. Their numerical magnitudes consequently cannot be compared directly with the intrinsic round-metric residuals in \cref{sec:results}.

\subsection{An explicit period model for \texorpdfstring{$S^4\times T^2$}{S4 x T2}}\label{sec:s4t2-periods}

The period model we study begins with complex two-tori over $t\in\C^*$ whose period columns, in additive coordinates $(z,w)$, are
\begin{equation}\label{eq:s4t2-periods-new}
 e_1=(2\pi,0),\quad e_2=(-it,1-t^2),\quad
 e_3=(-t,i(1+t^2)),\quad e_4=(-i\log t,2t).
\end{equation}
The logarithm is taken locally; continuation around the origin sends $e_4$ to $e_4+e_1$. The displayed real determinant is
\begin{equation}
 \det\nolimits_\R(e_1,e_2,e_3,e_4)
 =2\pi(1+r^2)\bigl[2r^2-(1-r^2)\log r\bigr],\qquad r=|t|,
\end{equation}
which is positive for every $r>0$. Thus the regular family is non-degenerate, and holomorphy gives local integrability. The model also specifies a chart at infinity, two cusp fillings and a logarithmic transform. Their topological consequences, including the identification with $S^4\times T^2$, are not established by these residual calculations.

A further experiment varies the periods within the prescribed polynomial form. Writing $\zeta=e^{iz}$, the translations have $w$-components $u_j=a_j+c_jt^2$ and multiplicative exponents $\alpha_j=-c_jt$ for $j=2,3$, while $u_4=b_4t$. Together with a constant $\beta$, these give twelve real parameters. The form satisfies the prescribed holomorphy and vanishing conditions at the two cusps. With $|a_2|^2+|a_3|^2+|c_2|^2+|c_3|^2=4$, the solver enforces opposite signs of the two cusp area forms and positivity of the normalised determinant $\det_\R(e_1,e_2,e_3,e_4)/[2\pi(1+r^4)(1+|\log r|)]$. These constraints retain the geometric ansatz while leaving the coefficients to be determined.

\begin{figure}[t]
 \centering
 \includegraphics[width=\textwidth]{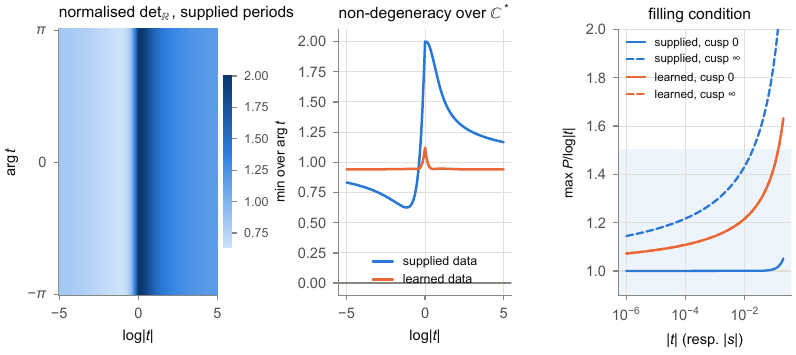}
 \caption{Determinant and cusp diagnostics for the supplied $S^4\times T^2$ period model and the selected learned candidate. The middle panel reports minima over the tested angles; the right panel evaluates the filling ratio near each cusp. Its scalar height function $P$ is unrelated to the tangent projector used in the ambient formulation. The displayed positive margins are numerical observations over the tested ranges, with the compactification and its topology requiring separate analytic arguments.}
 \label{fig:s4t2-periods-new}
\end{figure}

Nineteen of twenty random starts pass the recorded feasibility tests. All nineteen reach the largest feasible margin on the tested schedule, approximately $0.9411$; the next scheduled value, approximately $1.0822$, is infeasible. This is not a demonstrated optimum: the normalisation gives an upper bound of $1$ for a common margin on the two cusp area forms. The best learned candidate has smallest detected normalised determinant approximately $0.941$, compared with $0.626$ for \eqref{eq:s4t2-periods-new}. \Cref{fig:s4t2-periods-new} shows this comparison. The left panel maps the normalised determinant of the supplied periods over $\C^*$, and the middle panel compares its minimum over $\arg t$ for the two sets of periods. The learned candidate is nearly flat in $\log|t|$, while the supplied periods dip to about $0.63$. These results indicate that numerical non-degeneracy is not confined to a single isolated parameter choice. They do not establish positivity over all of $\C^*$ from finite samples, nor do they establish the compactified topology for the general learned parameters.

The recorded monodromy, infinity-chart and logarithmic-transform discrepancies are small for both sets of periods, at most $7.1\times10^{-12}$, with regular-part Nijenhuis residuals at most $1.1\times10^{-13}$ (\cref{tab:reconstruction-new}, lower block). The collar checks use $u=\log t$, and hence exclude the central fibre. Neither the extension through the fillings nor the identification with $S^4\times T^2$ is established here, for either set of periods; the right panel of \cref{fig:s4t2-periods-new} evaluates the filling ratio near each cusp for both sets of periods, and these are the quantities entering the fillings, which the tests themselves do not verify. The numerical outcome is a family of feasible regular period data, together with transition diagnostics that motivate further analysis. As discussed in \cref{sec:discussion}, this success demonstrates how much the representation can simplify the computation for the numerical reconstruction.

\section{Discussion}\label{sec:discussion}
The experiments show that learning almost complex structures and learning integrability have different numerical demands. The ambient representation resolves compatibility of local tensor descriptions by construction, leaving a global algebraic problem that can be checked independently of the differential equation. The almost complex investigations behave consistently across the selected controls: all nine admissible examples pass the numerical tests for both seeds, while $S^4$ does not and retains an order-one local defect. The large margin between the accepted residuals and the correction threshold makes these results practically useful, although of course the numerical search for a maximum remains distinct from a uniform proof.

For integrability, the strongest result is the performance of a small polynomial model with a second-order solver on the classical examples. Quadratic ambient matrices can express several familiar structures, but that observation alone would not ensure that optimisation from a learned starting field finds an integrable candidate. The repeated convergence on both seeds, followed by correction and fresh-point validation, supplies that additional numerical evidence. The comparison also cautions against identifying the size of a network with its effectiveness. A larger conjugator may face a more complicated representation problem and a less favourable optimisation problem than a directly parametrised field.

The difficult six-dimensional searches show the limits of this success. In particular, a plateau in the learning loss residuals indicating staled convergence may be caused by a restricted function class or alternatively by limited computational resources and an unfavourable initial condition. These experiments hence cannot be used to conclude approximation lower bounds for the polynomial or neural classes, or non-existence results. The recent fibration constructions make this last point especially concrete: geometric existence and numerical accessibility in a particular set of coordinates are separate questions.

Defect concentration gives a more specific explanation for some of the observed failures. Suppose, as a model calculation, that a bounded family of almost complex fields is integrable outside a tube of radius $\varepsilon$ around a subset of codimension $c$, and that its transition across the tube has derivative of size $O(\varepsilon^{-1})$. Since the Nijenhuis tensor is first order, this gives the estimate
\begin{equation}\label{eq:concentration-scaling}
\int_M\norm{N_{J_\varepsilon}}^2\,d\mathrm{vol}
=O(\varepsilon^{c-2}).
\end{equation}
Thus an average squared residual can decrease when $c>2$ without a comparable improvement in pointwise control. The estimate assumes the indicated family exists; it is not, by itself, a construction of smooth fillings on every manifold. It nonetheless explains why volume-weighted training can favour the observed concentration. The fact that this also occurs on the classical Hopf example rules out interpreting concentration as a signature of non-existence.

The reconstruction experiments answer a different and complementary question. Once periods are represented by holomorphic functions, local integrability is built in, and optimisation can focus on the much smaller compatibility problem. The very small reconstruction residuals achieved therefore reflect the benefit of favourable geometric parametrisations, and indicate strong numerical corroboration of existence of complex structures on these manifolds. 

Several extensions follow directly from these limitations. A rigorous treatment of accepted ACS candidates would combine their substantial numerical margins with computable global residual bounds. Integrability would require either a validated correction theorem for an approximate solution, or an exact description recovered from the learned field. A wider seed study and controlled comparisons of optimiser and representation would clarify the source of the classical successes. For the difficult cases, adaptive sampling should be paired with representations adapted to singular fibres and transition maps, while monitoring conditioning and algebraic admissibility. The present results establish a useful numerical framework for those investigations and identify precisely where further geometric information is needed.
Beyond this, this work introduces a series of numerical methods which have been validated for general search of almost complex and complex structures; applying this methodology to manifolds where existence of such structures is unknown remains the ultimate and most exciting development direction.

\appendix
\section{The orientation index}\label{app:orientation}
The orientation penalty in \cref{eq:orientation} is based on a pointwise fact that does not assume orthogonality. Let $A^2=-I$ on an oriented Euclidean vector space of dimension $2n$ and set $H=I+A^\top A$. Then $H$ defines a positive definite $A$-invariant inner product $h(u,v) = u^\top H v$, and
\begin{equation}
A^\top H=A^\top-A.
\end{equation}
The alternating form $\omega(u,v)=u^\top(A^\top-A)v$ is therefore $h(Au,v)$. In particular, $\omega(u,Au)=h(Au,Au)>0$ for non-zero $u$, so its symplectic orientation is the complex orientation of $A$.

To bound the Pfaffian, use $\det A=1$ and obtain
\begin{equation}
\det(A^\top-A)=\det H=\prod_{j=1}^{2n}(1+\sigma_j^2),
\end{equation}
where $\sigma_j$ are the singular values of $A$. Because $A^{-1}=-A$, the multiset of singular values is invariant under inversion. Pairing $\sigma$ with $\sigma^{-1}$ gives
\begin{equation}
(1+\sigma^2)(1+\sigma^{-2})=2+\sigma^2+\sigma^{-2}\ge4.
\end{equation}
It follows that $|\Pf(A^\top-A)|\ge2^n$. This proves both the sign and magnitude statements used in the penalty. A positively oriented change of orthonormal tangent frame leaves the value unchanged, so a single global tangent frame is unnecessary. For a general candidate $K$, the index can pass through zero; the penalty discourages this during optimisation but does not replace the algebraic residual.

\section{Implementation settings and provenance}\label{app:implementation}
\subsection{Main settings}
\Cref{tab:settings} records the principal settings of the campaign. The two main random seeds are 7 and 8. An additional $S^3\times S^1$ ACS seed (10) supplies a second base in the Hopf class for \cref{tab:supervised}, whose $S^5\times S^1$ fits use network seeds 7, 8, 17 and 18. Comparisons on the difficult cases are exploratory and do not all have matched budgets, so their timings and final residuals should be interpreted as results of the stated runs rather than a controlled ranking of optimisers.

\begin{table}[htbp]
\centering\small
\caption{Principal settings of the direct searches. Learning-rate arrows indicate cosine decay. The final admissibility checks are numerical, with the interpretation given in \cref{sec:validation}.}
\label{tab:settings}
\renewcommand{\arraystretch}{1.13}
\begin{tabularx}{\textwidth}{@{}l>{\raggedright\arraybackslash}X@{}}
\toprule
Stage & Settings\\
\midrule
Arithmetic and sampling & Double precision on CPU; independent round-product samples; usual Adam batch size $192$.\\
ACS polynomial attempt & $3000$ Adam steps; learning rate $5\times10^{-3}\to2\times10^{-4}$; coefficient standard deviation $0.3/\sqrt{n_{\mathrm{feat}}}$.\\
ACS retries & Within each class: plain objective, then Dirichlet coefficient $3\to10^{-3}$, then orientation penalty of weight $1$, each from a fresh initialisation; preliminary ascent uses $64$ starts and $50$ steps.\\
Neural refinement & GELU; hidden widths $256,256,256,256$ for $d=6$ and $128,128,128$ otherwise; output scale $0.02$; $10\,000$ Adam steps, $10^{-3}\to10^{-5}$; gradient clipping at $1$; kept only if it lowers the preliminary ascent maximum of the accepted polynomial field.\\
Final ACS refinement & Up to $200$ L-BFGS iterations on $768$ points, rejected if held-out loss increases; $30$ Newton--Schulz steps for correction.\\
Final algebraic check & $\radv(K)<1$ from $256$ starts and $150$ steps, followed by $\radv(J^*)<10^{-8}$ from $128$ starts and $75$ steps.\\
Neural conjugation & One MLP with two matrix outputs; output scale $0.5$ (exponential), $0.05$ (affine), $0.1$ (tail-aware); $1500$--$10\,000$ Adam steps, $10^{-3}\to10^{-5}$ ($3\times10^{-4}$ for the $S^4\times T^2$ continuation); L-BFGS chunks of $100$ iterations on $1024$ points, kept only if the held-out loss decreases.\\
Polynomial CS & Degree two, then degree three when needed; own ACS projected by least squares; $2000$ initial Adam steps, $3\times10^{-3}\to10^{-4}$.\\
LM rounds & Eight iterations per round; at least $128$ points and two residual entries per coefficient; fresh collocation each round and $2000$ held-out points; damping divided by $3$ on success and multiplied by $4$ on failure, carried between rounds within $[10^{-8},10^{6}]$.\\
CS stopping and checks & A degree is accepted when the held-out loss is below $10^{-16}$ and $\Nadv<10^{-6}$ ($64$ starts, $60$ steps); LM rounds stop early below $10^{-20}$; final dense summaries on $6000$ fresh points.\\
\bottomrule
\end{tabularx}
\end{table}

For the translation-invariant $S^4\times T^2$ search, one works over $x\in S^4\subset\R^5$ with the ambient projector $\diag(I_5-xx^\top,I_2)$. Coefficients depend only on $x$, while the two flat directions represent the commuting unit circle fields. A degree-$q$ polynomial matrix in this reduced model has $49\binom{5+q}{q}$ coefficients. The counts for degrees two, three and four are $1029$, $2744$ and $6174$, respectively. This accounts for the smaller parameter counts without imposing an integrable target.

\subsection{Records, checks and limits of reproducibility}
All numerical results in this paper come from one campaign run with the code in the repository. The driver \path{experiments/run_all.py} executes 72 jobs; four further $S^3\times S^1$ jobs (two additional ACS seeds, the homotopy-class computation and one supervised fit) were run with the same code. Every table, figure and quoted number is generated from the per-run records by \path{experiments/make_figures.py} and \path{experiments/paper_numbers.py}, and the repository's test suite (63 tests) passes on this code.

The records include useful consistency checks of the geometric operators. The classical octonionic ACS on the unit six-sphere has $\norm{N_J}^2=384$ in the ordered-pair convention of \cref{eq:parseval}. Explicit Hopf and Calabi--Eckmann fields test both the algebraic and differential residuals, while the formulas \cref{eq:hopfmatrix,eq:cematrix} test the expressive power of the degree-two model. The reduced invariant formulation is also compared with the full ambient product representation. These checks provide independent diagnostics of conventions and implementation; they are separate from the target-free optimisation.

Recorded timings are elapsed times of individual jobs, often sharing the same machine with other jobs or unrelated work. Summing them does not give the total campaign duration or a count of CPU-hours. For example, the jobs sum to approximately $6800$ job-minutes, of which approximately $5200$ are the supervised fits of \cref{tab:supervised}. The per-run timings in \cref{tab:classical-learning} are sufficient for the limited comparison made here, with the preceding ACS cost explicitly excluded.

\bibliographystyle{utphys}
\bibliography{references}

\end{document}